\documentclass[11pt]{article}

\usepackage{amssymb,amsfonts,amsmath,nicefrac,color,graphicx,setspace}

\graphicspath {{figures/}}

\usepackage{algorithm,algorithmic,enumitem}

\def\sT {{\sf T}}

\newcommand{\cA}{{\mathfrak A}}

\newcommand{\cM}{{\mathfrak M}}

\newcommand{\cq}{{\mathfrak q}}

\newcommand{\cV}{{\mathfrak V}}
\newcommand{\uV}{\underline{ V}}

\newtheorem{remark}{Remark}[section]

\newcommand{\V}{{\cal V}}
\newcommand{\Z}{{\cal Z}}
\newcommand{\F}{{\cal F}}

\newcommand{\X}{{\cal X}}

\newcommand{\be}{\begin{equation}}
\newcommand{\ee}{\end{equation}}

\newcommand{\rr}{\rightrightarrows}

\def\w{\omega}

\def\O{\Omega}

\newcommand{\avr}{{\sf AV@R}}

\def\valrisk {{\sf V@R}}

\def\bbr{{\Bbb{R}}} 
\def\bbe{{\Bbb{E}}} 

\newcommand{\ip}[1]{\langle #1 \rangle}
\newcommand\tab[1][1cm]{\hspace*{#1}}

\newcommand{\ind}{{\mbox{\boldmath$1$}}}

\begin{document}
\title{\textbf{Risk Neutral Reformulation Approach to Risk Averse  Stochastic Programming}}

\author{
{\bf  Rui Peng Liu}   and
 {\bf Alexander  Shapiro}\thanks{
Research of this author  was partly supported by NSF grant 1633196.}\\
School of Industrial and Systems Engineering\\
Georgia Institute of Technology\\
Atlanta, GA 30332-0205\\
}
\maketitle
\thispagestyle{empty}

\begin{abstract}
The aim of this paper is to show that in some cases risk averse multistage stochastic programming problems can be reformulated in a form of risk neutral setting. This is achieved by a change of the reference probability measure making ``bad" (extreme) scenarios more frequent. As a numerical example we demonstrate advantages of such change-of-measure approach applied to the Brazilian Interconnected Power System operation planning problem.
\end{abstract}


\newpage

\pagenumbering{arabic}

\setcounter{equation}{0}
\section{Introduction}\label{sec:intr}

There are many practical problems where one has to make decisions sequentially based on
data (observations) available at time of the decision. In the stochastic programming approach the underlying  data  is modeled as a random process with a specified probability distribution. We can refer to the books \cite{pfl14},\cite{SDR}, and references therein, for a thorough discussion of the Multistage Stochastic Programming (MSP). In the risk neutral formulation of MSP problems,   expected value of the  total cost is supposed to be  optimized (minimized). Of course for a particular realization of the data process the corresponding cost can be quite different from its average. This motivates to consider risk averse approaches where one tries to control   high costs by imposing some type  of penalty on high cost realizations of the data process.

An axiomatic approach to risk   was suggested in the pioneering paper by   Artzner et al  \cite{Artzner},   where   the concept of   coherent risk measures was introduced.  By the  Fenchel - Moreau theorem coherent risk measures can be represented in the dual form as maximum of expected values. Consequently  optimization problems involving coherent risk measures can be written in a minimax form. This suggests that such risk averse problems can be formulated as risk neutral problems with respect to an appropriate worst case  probability distribution (cf., \cite[Remarks 24-25]{SDR}).
Although such worst case probability distribution is not known a priori, we are  going to demonstrate that in some cases it can be approximated  in a computationally feasible way.

From several points of view a natural example of   coherent risk measures is the so-called Average Value at Risk ($\avr$) (under different  names, such as  Conditional Value at Risk,  Expected Shortfall,  Expected Tail Loss, this    was discovered   and rediscovered  in various equivalent forms  by several   authors over many years). A nested risk averse approach, using  convex combinations  of the expectation and $\avr$,  was suggested  in \cite{ejor}  for  controlling high costs  in  planning of hydropower generation. In that approach risk of high costs is controlled  by imposing an appropriate  penalization on such high costs    at every stage of the decision process conditional on observed realizations  of the random data process.
One of the criticisms   of such risk averse approach  is that
the corresponding objective is formulated  in a nested form and is
difficult for  an intuitive interpretation.

The contribution of this paper is twofold. We demonstrate that in some situations it is possible to reformulate the considered risk averse problem in a risk neutral form by making an appropriate change of the probability  measure.
This leads to an intuitive interpretation of controlling the risk by giving higher  weights to ``bad scenarios".  The idea of constructing scenario trees with extreme (bad) scenarios was considered before (e.g., \cite[Chapter 2]{zie}). In that respect our approach is quite different.  We relate it to the modern risk averse approach to MSP and blend it with the Stochastic Dual Dynamic Programming (SDDP)  type algorithm. In particular this allows construction of  lower and {\em upper} numerical bounds, for the constructed risk neutral problem, following the standard risk neutral methodology of the SDDP method. Also our approach is different from  the approach based on extended polyhedral risk measures which allows reformulation of the respective  risk-averse problems as risk-neutral problems (with additional constraints and variables) (cf., \cite{guigrom:2012}).

Another contribution of the suggested ``change of the probability  measure" method is an improvement of the rate of convergence of the straightforward risk averse method. The idea is somewhat related to the classical Importance Sampling  techniques  although is not exactly the same.  The standard Importance Sampling methodology is aimed at reducing variance of the respective sample estimates, while we are concerned with an approximation of the corresponding risk averse problem. An intuitive explanation of our approach is that by generating ``bad" scenarios more often one becomes more conservative and risk aware  in his/her decisions.

This paper is organized as follows. In the next section we discuss a static case of the risk averse stochastic programming. We show how in some situations  the corresponding worst case distribution can be computed. In section \ref{sec:mult} we extend this to a multistage setting.   In section \ref{sec:sddp} we discuss a  risk averse variant of the Stochastic Dual Dynamic Programming (SDDP)  algorithm and its reformulation in a risk neutral form. In section \ref{sec:ne} we give a numerical example based on the Brazilian Interconnected Power System operation planning problem. Finally   section \ref{sec:conc} is  devoted to  concluding remarks.

\setcounter{equation}{0}
\section{Static case}
\label{sec:stat}
Consider the following risk averse stochastic program
\begin{equation}\label{sta-1}
 \min_{x\in \X}\rho[F(x,\w)],
\end{equation}
where $(\O,\F,P)$ is a probability space,  $\X\subset \bbr^n$,   $F:\bbr^n\times \O\to \bbr$, and $\rho:\Z\to\bbr$ is a coherent risk measure defined on a linear space $\Z$ of random variables $Z:\O\to\bbr$.
We assume that for every $x\in \X$, random variable $F_x(\w)=F(x,\w)$ belongs to $\Z$.
We also assume that problem \eqref{sta-1} is convex, i.e., the set $\X$ is convex and $F(x,\w)$ is convex in $x$ for a.e. $\w\in \O$.
In particular we deal with risk measures of the form
\begin{equation}
\label{sta-2}
 \rho(Z):=(1-\lambda)\bbe[Z]+\lambda\avr_\alpha (Z),\;\lambda \in   (0,1),
\end{equation}
where  $\Z=L_1(\O,\F,P)$ paired with its dual space $\Z^*=L_\infty(\O,\F,P)$,
and
\[
\avr_\alpha (Z)=\inf_{u\in \bbr}\left\{u+\alpha^{-1}\bbe[Z-u]_+\right\}, \;\alpha\in (0,1).
\]
In the above variational form,  $\avr$ was defined in \cite{rock2000} under the name ``Conditional Value at Risk".

By the    Fenchel - Moreau theorem, real valued
coherent risk measure $\rho$ has dual representation (cf., \cite{RuSh:2006a})
 \begin{equation}
\label{sta-3}
 \rho (Z)=\sup_{\zeta\in \cA}\bbe_\zeta [Z],
\end{equation}
where $\cA\subset \Z^*$ is a convex weakly$^*$ compact  set of probability  density functions   and
\[
\bbe_\zeta [Z]:=\int_{\O}     Z(\w) \zeta(\w) dP(\w),
\]
is the expectation with respect to the probability  measure $dQ=\zeta dP$. Hence problem \eqref{sta-1} can be written in the following minimax form
\begin{equation}\label{sta-4}
 \min_{x\in \X}\sup_{\zeta\in \cA}\bbe_\zeta[F_x].
\end{equation}
A dual of problem \eqref{sta-4} is obtained by interchanging the `min' and `max' operators:
\begin{equation}\label{sta-5}
 \max_{\zeta\in \cA}\inf_{x\in \X}\bbe_\zeta[F_x].
\end{equation}

A point $(\bar{x},\bar{\zeta})\in \X\times\cA$ is said to be a saddle point of the above minimax problems if
\begin{equation}
\label{sta-5a}
 \bbe_{\bar{\zeta}}[F_x]\ge \bbe_{\bar{\zeta}}[F_{\bar{x}}]\ge  \bbe_{\zeta}[F_{\bar{x}}],\;\forall (x,\zeta)\in \X\times   \cA.
 \end{equation}
Under mild regularity conditions the minimax problem \eqref{sta-4} has a saddle point $(\bar{x},\bar{\zeta})\in \X\times\cA$. Then $\bar{x}$ is an optimal solution of problem \eqref{sta-4},  $\bar{\zeta}$ is an optimal solution of problem \eqref{sta-5}, optimal values of problems \eqref{sta-4} and \eqref{sta-5} are equal to each other and are equal to the optimal value of the following problem
\begin{equation}\label{sta-6}
\inf_{x\in \X}\bbe_{\bar{\zeta}}[F_x].
\end{equation}
It follows from the first inequality in \eqref{sta-5a} that if $\bar{x}$ is an optimal solution of problem \eqref{sta-1} (i.e.,
$(\bar{x},\bar{\zeta})$ is  a saddle point), then $\bar{x}$ is also  an optimal solution of problem \eqref{sta-6}. That is, the set of optimal solutions of problem \eqref{sta-6} contains the set of optimal solutions of problem \eqref{sta-1} (it can happen that  the set of optimal solutions of problem \eqref{sta-6} is larger than  the set of optimal solutions of problem \eqref{sta-1}).

That is,  risk averse problem \eqref{sta-1} can be formulated as risk neutral problem \eqref{sta-6} with respect to the ``worst"  probability measure $dQ=\bar{\zeta}dP$. Of course  $\bar{\zeta}$ is not known, its evaluation requires solution of the minimax problem \eqref{sta-5}. Nevertheless this gives us a direction for constructing approximation of problem \eqref{sta-6}. For $Z=F_{\bar{x}}$
we have that
\begin{equation}\label{sta-7}
\bar{\zeta}\in \arg\max_{\zeta\in\cA} \bbe_{\zeta}[Z].
\end{equation}
Recall that for $Z\in \Z$,
\begin{equation}\label{sta-8}
 \arg\max_{\zeta\in\cA} \bbe_{\zeta}[Z]=\partial \rho(Z)
\end{equation}
(cf. \cite[eq.(6.43), p.265]{SDR}).

Consider   risk measure \eqref{sta-2}.  We have that
\begin{equation}\label{sta-9}
\partial \rho(Z)=(1-\lambda)\{\ind\}+\lambda\partial (\avr_\alpha)(Z)
\end{equation}
with $\ind(\cdot)\equiv 1$ and  (cf. \cite[eq.(6.74), p.273]{SDR})
\begin{equation}\label{sta-10}
\partial (\avr_\alpha)(Z)=\left\{\zeta:\bbe[\zeta]=1,
\begin{array}{lll}
 \zeta(\w)=\alpha^{-1},&\text{if}& Z(\w)>\valrisk_{\alpha}(Z), \\
\zeta(\w)=0,&\text{if}& Z(\w)<\valrisk_{\alpha}(Z), \\
\zeta(\w)\in[0,\alpha^{-1}],&\text{if}& Z(\w)=\valrisk_{\alpha}(Z),
\end{array}
\right.
\end{equation}
where
\[
\valrisk_{\alpha}(Z)=\inf\{t:P(Z\le t)\ge 1-\alpha\}.
\]

Suppose further   that the space  $\O=\{\w_1,...,\w_N\}$  is finite equipped with equal probabilities $p_i=1/N$, $i=1,...,N$.  Then denoting $Z_i=Z(\w_i)$,  random variable $Z:\O\to \bbr$ can be identified with vector $(Z_1,...,Z_N)\in \bbr^N$.  In that case
\begin{equation}\label{avr1}
 \rho(Z)=\frac{(1-\lambda)}{N} \sum_{i=1}^N Z_i+
\lambda \bar{u}+\frac{\lambda}{\alpha N}\sum_{i=1}^N [Z_i-\bar{u}]_+,
\end{equation}
where $\bar{u}=\valrisk_{\alpha}(Z)$.  Let $Z_{(1)}\le...\le Z_{(N)}$ be values $Z_i$, $i=1,...,N$,  arranged in the increasing order. Then $\bar{u}=Z_{(\kappa)}$, where
$\kappa:=\lceil (1-\alpha)N \rceil$ with  $\lceil a \rceil$ denoting   the smallest integer $\ge a$. We assume that $N$ is large enough so that  $\kappa\le N-1$. Then we can write
\begin{equation}\label{avr2}
 \rho(Z)=\frac{(1-\lambda)}{N} \sum_{i=1}^{N} Z_{(i)}+\lambda  Z_{(\kappa)}+ \frac{\lambda}{\alpha N}\sum_{i=\kappa+1}^N(Z_{(i)}-Z_{(\kappa)})=
 \sum_{i=1}^N \cq_i Z_{(i)},
\end{equation}
where
\begin{equation}\label{distr}
\cq_i:=\left\{
\begin{array}{lll}
(1-\lambda)/N &{\rm if }&i<\kappa, \\
(1-\lambda)/N+\lambda -\lambda(N-\kappa)/(\alpha N) &{\rm if }& i=\kappa, \\
(1-\lambda)/N+\lambda/(\alpha N) &{\rm if } & i>\kappa.
\end{array}\right .
\end{equation}
Note that $\cq_i\ge 0$, $i=1,...,N$,  and $\sum_{i=1}^N \cq_i=1$.
That is, in order to find the worst  probability measure $Q$  we only  need to identify the ``bad" outcomes of the distribution of $F_{\bar{x}}$, i.e., which
  values $Z_i=F_{\bar{x}}(\w_i)$  are larger than $\valrisk_{\alpha}(F_{\bar{x}})$.
We sometimes write $\cq_i=\cq_i(Z)$ for values defined in \eqref{distr} associated with vector $Z=(Z_1,...,Z_N)$.

\setcounter{equation}{0}
\section{Multistage  case}
\label{sec:mult}

Consider a multistage  risk averse  stochastic
programming problem given   in the
following nested form (cf.,  \cite{SDR})
\begin{equation}
\label{cm-1}
\min_{x_1\in \X_1}  f_1(x_{1}) +
\rho_{2|\xi_{[1]}}\left [ \inf_{x_2\in \X_2(x_1,\xi_2) }
 f_2(
x_{2},\xi_2) +  \cdots +\rho_{T|\xi_{[T-1]}} \,\Big[
\inf_{x_T\in \X_T(x_{T-1},\xi_T)} f_T(x_{T},\xi_T)
 \Big]\right],
\end{equation}
driven by the random data  process
$\xi_1,\xi_2,{\dots},\xi_T$.  Here $x_t\in
\bbr^{n_t}$, $t=1,{\dots},T$, are  decision
variables, $f_t:\bbr^{n_{t}}\times \bbr^{d_t}\to
\bbr$ are continuous functions,
$\X_t:\bbr^{n_{t-1}}\times \bbr^{d_t}\rr
\bbr^{n_{t}}$, $t=2,{\dots},T$,  are measurable
closed valued  multifunctions and $\rho_{t|\xi_{[t-1]}}$ are conditional coherent  risk mappings  (we use notation $\xi_{[t]}:=(\xi_1,...,\xi_t)$ for the history of the data process).
As the main example we consider the following  conditional counterpart   of the risk measure \eqref{sta-2}:
\begin{equation}
\label{cm-crm}
\rho_{t|\xi_{[t-1]}}[Z]=(1-\lambda) \bbe\left [Z\big |\xi_{[t-1]}\right]+ \lambda \avr_\alpha \left [Z\big |\xi_{[t-1]}\right],\;\lambda,\alpha\in (0,1).
\end{equation}
The first stage
data, i.e.,  the vector $\xi_1$,  the function
$f_1:\bbr^{n_1}\to\bbr$, and the
  set  $\X_1\subset \bbr^{n_1}$ are deterministic.

We assume that  problem \eqref{cm-1} is convex, i.e., functions $f_t(\cdot,\xi_t)$    and sets $\X_t(x_{t-1},\xi_t)$ are convex.
 It is said that  the multistage problem \eqref{cm-1} is {\em linear} if the objective functions and the
constraint functions are linear, that is
\begin{equation}
\label{cm-2}
\begin{gathered}
f_t(x_t,\xi_t):=  c_t^{\sT} x_t,\;\; \X_1:=\left \{x_1:A_{1}x_1= b_1,\;  x_1  \ge 0\right \},\\
\X_t(x_{t-1},\xi_t):=\left \{x_t:B_{t}x_{t-1}+A_{t}x_t=  b_t,\;   x_t  \ge 0\right \},\;t=2,{\dots},T.
 \end{gathered}
\end{equation}
Here, $\xi_1:=(c_1,A_{1},b_1)$  is known at the first
stage (and hence is nonrandom), and
$\xi_t:=(c_t,B_{t},A_{t},b_t)\in \bbr^{d_t}$,
$t=2,{\dots},T$,    are data vectors some
(or all) elements of which can be random. Linear problems are convex.

Problem \eqref{cm-1} can be written in the following equivalent form
\begin{equation}
\label{cm-3}
\min_{\pi\in \Pi}
\bar{\rho}\left[  f_1(x_{1}) +f_2(x_{2}(\xi_{[2]}),\xi_2) +  \cdots +
  f_T(x_{T}(\xi_{[T]}),\xi_T)\right],
\end{equation}
where $\Pi$ is the set of policies $\pi=\left(x_1,x_{2}(\xi_{[2]}),...,x_{T}(\xi_{[T]})\right)$ satisfying the feasibility constraints of problem \eqref{cm-1},  and $\bar{\rho}$ is the composite risk measure (cf., \cite[p.318]{SDR})
\begin{equation}
\label{nestf}
\bar{\rho}[Z]=\rho_{2|\xi_{[1]}}\Big[\rho_{3|\xi_{[2]}}\big[\cdots \rho_{T|\xi_{[T-1]}}[Z]\big]\Big].
\end{equation}
 The composite risk measure $\bar{\rho}$ is given in the nested form \eqref{nestf}  and   could be quite  complicated. As it was already mentioned in the introduction,  this may raise an objection of an intuitive interpretation of the overall objective of the risk averse formulation \eqref{cm-3}.

 Anyway risk measure $\bar{\rho}$ is coherent
 and has the dual representation
 \begin{equation}
\label{cm-4}
\bar{\rho}[Z]=\sup_{Q\in \cM}\bbe_Q[Z],
\end{equation}
where $\cM$ is a set of probability measures (distributions) of $\xi_{[T]}$ {\em absolutely continuous} with respect to the reference probability measure of the data process.
Therefore problem \eqref{cm-3} can be written in the following minimax form
\begin{equation}
\label{cm-5}
\min_{\pi\in \Pi}\max_{Q\in \cM}\bbe_Q
\left[  f_1(x_{1}) +f_2(x_{2}(\xi_{[2]}),\xi_2) +  \cdots +
  f_T(x_{T}(\xi_{[T]}),\xi_T)\right].
\end{equation}

The dual of problem \eqref{cm-5} is the problem
\begin{equation}
\label{cm-6}
\max_{Q\in \cM}\min_{\pi\in \Pi}\bbe_Q
\left[  f_1(x_{1}) +f_2(x_{2}(\xi_{[2]}),\xi_2) +  \cdots +
  f_T(x_{T}(\xi_{[T]}),\xi_T)\right].
\end{equation}
Under mild regularity conditions optimal values of problems \eqref{cm-5} and \eqref{cm-6} are equal to each other. Suppose further that problem \eqref{cm-6} has an optimal solution $\bar{Q}$. Then problem \eqref{cm-1} has the same optimal value as the risk neutral  problem
\begin{equation}
\label{cm-7}
 \min_{\pi\in \Pi}\bbe_{\bar{Q}}
\left[  f_1(x_{1}) +f_2(x_{2}(\xi_{[2]}),\xi_2) +  \cdots +
  f_T(x_{T}(\xi_{[T]}),\xi_T)\right],
\end{equation}
and the set of optimal solutions of problem
\eqref{cm-1} is contained in the set of optimal solutions of problem \eqref{cm-7}.
When the number of scenarios is finite, i.e., the data process can be represented by a finite scenario tree,   the change of measure can be described in a constructive way (cf., \cite[Remarks 24-25, pp.314-315] {SDR}).

We assume in the remainder of this section  that the data process $\xi_t$ is {\em stagewise independent}, i.e., random vector $\xi_{t+1}$ is independent  of $\xi_{[t]}$, $t=1,...,T-1$ (although  $\xi_1$ is supposed to be deterministic, we include it for uniformity of the notation.). Suppose further that   the  conditional  risk mappings $\rho_{t|\xi_{[t-1]}}$ are given as the conditional counterparts  of coherent risk measures $\rho_t$.
 In the stagewise independent case the joint probability distribution of $(\xi_1,...,\xi_T)$ is determined by the marginal distributions of each
$\xi_t$, $t=1,...,T$, and the corresponding cost-to-go (value) functions can be written as
\begin{equation}
\label{dim1}
V_{t}\left (x_{t-1},\xi_{t}\right )=\inf_{x_t\in
\X_t(x_{t-1},\xi_t)}\big\{f_t(x_t,\xi_t)+\V_{t+1}\left
(x_t\right )\big\},
\end{equation}
where
\begin{equation}
\label{dim2}
\V_{t+1}(x_t):=\rho_{t+1}\big(V_{t+1}
(x_t,\xi_{t+1} )\big),
\end{equation}
$t=1,...,T$, with $\V_{T+1}(\cdot)\equiv 0$.

Let
\begin{equation}
\label{dim3}
\bar{x}_t\in \arg\min_{x_t\in
\X_t(\bar{x}_{t-1},\xi_t)}\big\{f_t(x_t,\xi_t)+\V_{t+1}\left
(x_t\right )\big\}.
\end{equation}
Note that  $\bar{x}_t$,  $t=2,...,T$, is a function of   $\bar{x}_{t-1}$ and  $\xi_t$, and that the policy  $(\bar{x}_1,...,\bar{x}_T)$ is an optimal solution of the corresponding multistage problem \eqref{cm-1}.
For each $\rho_{t+1}$  we can consider its dual representation of the form \eqref{sta-3}
\begin{equation}
\label{dual}
 \rho_{t+1} (Z)=\sup_{\zeta_{t+1}\in \cA_{t+1}}\bbe_{\zeta_{t+1}} [Z],
\end{equation}
with the corresponding set $\cA_{t+1}$ of density functions. Consider  a saddle point   of the   minimax problem
\begin{equation}
\label{dim4}
 \min_{x_t\in
\X_t(\bar{x}_{t-1},\xi_t))}\sup_{\zeta_{t+1}\in \cA_{t+1}} \left\{f_t(x_t,\xi_t)+ \bbe_{\zeta_{t+1}} [V_{t+1}
(x_t,\xi_{t+1} )]  \right\}.
\end{equation}
For $x_t=\bar{x}_t$ we need to solve the problem
\begin{equation}
\label{dim5}
 \max_{\zeta_{t+1}\in \cA_{t+1}}   \bbe_{\zeta_{t+1}} [V_{t+1}
(\bar{x}_t,\xi_{t+1} )] ,
\end{equation}
 in order to find the  component $\bar{\zeta}_{t+1}$ of the corresponding  saddle point
$(\bar{x}_t,\bar{\zeta}_{t+1})$.

Suppose now that marginal distribution of $\xi_t$, $t=2,...,T$,  is discretized by generating $N$ points
$\xi_t^1,...,\xi_t^N$,  each assigned with the same  probability $1/N$. Suppose further that $\rho_t=\rho$, for all $t$, with $\rho$ given in the form \eqref{sta-2}.  Then as it was discussed in section \ref{sec:stat},  in order to find the corresponding worst case density (worst case distribution)  we only  need to identify the ``bad" outcomes of the distribution of $V_{t+1}(\bar{x}_t,\xi_{t+1} )$ (see \eqref{distr}),
 i.e., which values $V_{t+1}(\bar{x}_t,\xi^j_{t+1} )$ are larger than $\valrisk_{\alpha}(V_{t+1}(\bar{x}_t,\xi_{t+1} ))$.  Recall that $\bar{x}_t$,  $t=2,...,T$, is a function of   $\bar{x}_{t-1}$ and  $\xi_t$.  Nevertheless if   $\valrisk_{\alpha}(V_{t+1}(\bar{x}_t,\xi_{t+1} ))$ is stable with respect to different realizations of  $\bar{x}_{t}$, then it would be possible to assign    weights (probabilities) to $\xi_t^j$, at every stage $t$ of the decision  process, in such a way  that the constructed  risk neutral problem will be equivalent to the original risk averse multistage problem.  We will discuss such an example in the following sections.

\section{Stochastic Dual Dynamic Programming algorithm}\label{sec:sddp}

In the risk neutral setting  the Stochastic Dual Dynamic Programming (SDDP) algorithm was introduced in Pereira and   Pinto \cite{per91}.  Its origins can be traced to the nested decomposition algorithm of Birge \cite{bir85}. A risk averse variant of the SDDP method, based on a convex combination of the expectation and $\avr$,  was introduced  in \cite{sha:11}  and implemented in \cite{phil12}. For  extended polyhedral risk measures a variant of the SDDP method is discussed in \cite{guigrom:2012}.
Convergence of the SDDP  type algorithms is discussed in \cite{Gir:2014},\cite{guig:2016},\cite{gui2017}, and references therein.
{\color{black} For a view of the SDDP method as a form of approximate dynamic programming we refer to \cite{powell}}.

Consider the multistage stochastic programming problem \eqref{cm-1}. We assume that  the problem is linear, i.e., the data are given in the form \eqref{cm-2}. Moreover we assume   that  the data process $\xi_t$ is stagewise independent and  has a finite number of scenarios. That is, the marginal distribution of $\xi_t$, $t=2,...,T$, has $N$ realizations
 $\xi_t^1,...,\xi_t^N$ each having  the same  probability $1/N$  (for the sake of simplicity we assume that the number $N$ of discretization points at each stage $t$ is the same).
The SDDP algorithm is a cutting plane type method  designed to solve such  multistage convex stochastic programming problems.
Being an iterative approach, the SDDP algorithm
   progressively refines lower piecewise linear approximations of $\V_{t+1}(\cdot)$. It has two major steps at each iteration.


\paragraph{Forward Step}

At iteration $m$, suppose for each $t = 1, \ldots, T-1$, we have a finite set $\mathcal{S}_{t+1}^m$ of affine minorants $s(\cdot)$ of $\V_{t+1}(\cdot)$, then $\mathfrak{V}_{t+1}^m(\cdot) := \max_{s\in \mathcal{S}_{t+1}^m} s(\cdot)$ is a lower piecewise linear approximation of $\V_{t+1}(\cdot)$, and
\begin{equation} \label{approx1}
\uV_t^m (\cdot, \xi_t) := \inf_{x_t\in \mathcal{X}_t(\cdot, \xi_t)} \{f_t(x_t, \xi_t) +  \cV_{t+1}^m(x_t)\}
\end{equation}
is a lower convex approximation of $V_t(\cdot, \xi_t)$ for each $\xi_t$.

In the forward step of the algorithm a random realization  $\xi_t^{j_t}$, $j_t\in \{1,....,N\}$,  from the distribution of $\xi_t$, $t=2,...,T$, is sampled (we refer to this as the {\em uniform sampling}). In the case that each realization $\xi_t^j$ of $\xi_t$ has the same probability $1/N$, each scenario (sample path) $\{\xi_t^{j_t}\}_{t = 2, \ldots, T}$ of the data process is sampled with probability    $1/N^{T-1}$. At stage $t$  the corresponding  trial point $x_t^m$ is defined to be the  optimal solution
\begin{equation} \label{approx2}
x_t^m \in\arg\min\limits_{x_t\in \mathcal{X}_t(x_{t-1}^m, \xi_t^{j_t})}\left  \{f_t(x_t, \xi_t^{j_t}) +  \mathfrak{V}_{t+1}^m(x_t)\right \}.
\end{equation}
Note that $f_1(x_1^m) + \mathfrak{V}_2^m(x_1^m)$ gives a lower bound of the total cost.

\begin{remark}
\label{rem-sub}
{\rm
When the problem is risk neutral (i.e., $\lambda = 0$), the measure $\bar{Q}$ in \eqref{cm-7} is given by a discrete measure that assigns probability $1/N^{T-1}$  to each scenario $\{\xi_t^{j_t}\}$. Viewing $\{x_t^m\}$ as functions of $\{\xi_t^{j_t}\}$, the quantity $\mathbb{E}_{\bar{Q}} [\sum_{t = 1}^T f_t(x_t^m (\xi_t), \xi_t) ]$ is then an upper bound of the (optimal) total cost, and $\sum_{t = 1}^T f_t(x_t^m, \xi_t^{j_t})$ is the corresponding unbiased estimator. By  generating a number of  scenarios  one  can  construct a confidence interval for the upper bound.
In the risk averse setting construction of the corresponding upper bound is more involved. A straightforward analogue of the statistical upper bound, used in the risk neutral case,  does not work in the risk averse setting  (cf., \cite{ding}). An approach to constructing (nonstochastic)   upper bounds,  which is  also applicable in the risk averse setting, was suggested in
\cite{PhMatFin}.
}
\end{remark}


\paragraph{Backward Step}
Given trial  points $\{x_t^m\}_{t = 1, \ldots, T - 1}$ at iteration $m$,  the approximations $\mathfrak{V}_t^m(\cdot)$ are refined  sequentially from $t = T$ down to $t=2$. By equation \eqref{avr2}, we have
\begin{equation} \label{approx3}
\begin{split}
\rho_t(\uV_t^m(x_{t-1}, \xi_t)) &= \frac{1 - \lambda}{N}  \sum_{j = 1}^N \uV_t^m(x_{t-1}, \xi_t^j) + \lambda \uV_t^m(x_{t-1}, \xi_t^{(\kappa)}) \\
			&\tab \tab + \frac{\lambda}{\alpha N}  \sum_{j = \kappa + 1}^N\left(\uV_t^m(x_{t-1}, \xi_t^{(j)}) - \uV_t^m(x_{t-1}, \xi_t^{(\kappa)})\right ),
\end{split}
\end{equation}
where $\uV_t^m(x_{t-1}, \xi_t^{(1)})\le\ldots\le \uV_t^m(x_{t-1}, \xi_t^{(N)})$ are values $\uV_t^m(x_{t-1}, \xi_t^j))$, $j = 1, \ldots, N$, arranged in the increasing order, and $\kappa := \lceil(1 - \alpha) N\rceil$.
The corresponding  affine  minorant $s_t^m(\cdot)$ of $\rho_t(\uV_t^{m+1}(\cdot, \xi_t))$ at $x_{t-1}^m$
is constructed by computing  a subgradient $g_t^j$ of $\uV_t^{m+1}(\cdot, \xi_t^j)$ at $x_{t-1}^m$  (cf., \cite[Section 4.1]{ejor}).
Observe that $s_t^m(\cdot)$ is also an affine minorant of $\mathcal{V}_t(\cdot)$, since $\rho_t(\uV_t^{m+1}(\cdot, \xi_t))$ is a lower convex approximation of $\mathcal{V}_t(\cdot)$. Next, define $\mathcal{S}_t^{m+1} := \mathcal{S}_t^m\cup \{s_t^m\}$. Then  $\mathfrak{V}_t^{m+1}(\cdot) = \max\{\mathfrak{V}_t^m(\cdot), s_t^m(\cdot)\}$ is a refined lower piecewise linear approximation of $\mathcal{V}_t(\cdot)$.

We now discuss how to generate the subgradients  $g_t^j$ when the multistage problem is linear (see \eqref{cm-2}). Suppose $\xi_t^j = (c_t^j, B_t^j, A_t^j, b_t^j)$, and $\mathfrak{V}_{t + 1}^{m + 1}(\cdot)$ has been generated, where $\mathfrak{V}_{T+1}^{m+1}(\cdot)\equiv 0$. By definition \eqref{approx1},
\begin{equation} \label{approx5}
\uV_t^{m+1}(x_{t-1}, \xi_t^j) = \inf_{x_t\ge 0} \left \{(c_t^j)^{\sT} x_t + \mathfrak{V}_{t+1}^{m+1}(x_t) : B_t^j x_{t-1}+ A_t^j x_t = b_t^j\right \}.
\end{equation}
A subgradient of $\uV_t^{m+1}(\cdot, \xi_t^j)$ at $x_{t-1}^m$ is given by $g_t^j = -(B_t^j)^\top \pi_t^j$, where $\pi_t^j$ is an optimal solution to the dual problem of
\begin{equation} \label{approx6}
\min_{x_t\ge 0}\left  \{(c_t^j)^{\sT} x_t + \mathfrak{V}_{t+1}^{m+1}(x_t) : B_t^j x_{t-1}^m + A_t^j x_t = b_t^j\right \}.
\end{equation}
Problem \eqref{approx6} can be reformulated as the linear programming problem
\begin{equation} \label{approx7}
\begin{split}
\min_{x_t\ge 0, \theta}\ \ \{(c_t^j)^{\sT} x_t + \theta : B_t^j x_{t-1}^m + A_t^j x_t = b_t^j, s(x_t)\le \theta, s\in \mathcal{S}_{t+1}^{m + 1}\}. \\
\end{split}
\end{equation}
The dual problem of \eqref{approx7} is also a linear programming problem, hence a dual optimal solution $\tilde{\pi}_t^j$ can be obtained efficiently, and $\pi_t^j$ is exactly the entries of $\tilde{\pi}_t^j$ corresponding to the constraints $B_t^j x_{t-1}^m + A_t^j x_t = b_t^j$.

\begin{remark}
\label{rem-sev}
{\rm
The  backward steps can be carried out with arbitrary (feasible)  points $x_t^m$,  $t=1, \ldots, T-1$. The trial points generated at the  forward steps are  one of the  possibilities.
As we shall see later, there could be more efficient ways to generate the trial points.
}
\end{remark}

\begin{remark}
\label{rem-rarssddp}
{\rm Depending on whether $\lambda > 0$ or $\lambda = 0$ we refer to the SDDP algorithm as risk averse or risk neutral,  respectively.

}
\end{remark}

\begin{algorithm}[H]
\label{alg-1}
\caption{the SDDP algorithm for linear problem (uniform sampling)}
\begin{algorithmic}
\STATE $\kappa\leftarrow \lceil(1 - \alpha)N)\rceil, m\leftarrow 1, \mathfrak{V}_t^1(\cdot)\leftarrow 0\ \forall t$
\WHILE{termination criterion not met}

\STATE (Forward Step)
\STATE sample $\{\xi_t^{j_t}\}_{t = 2, \ldots, T - 1}$ according to  the original  distribution of the data process $\xi_t$.
\STATE $x_1^m\leftarrow \text{argmin}_{x_1\ge 0} \{c_1^{\sT} x_1 + \mathfrak{V}_2^m(x_1) : A_1 x_1 = b_1\}$
\FOR{$t = 2, \ldots, T - 1$}
\STATE $x_t^m \leftarrow \text{argmin}_{x_t\ge 0} \{(c_t^{j_t})^{\sT} x_t + \mathfrak{V}_{t + 1}^m(x_t) : B_t^{j_t} x_{t - 1}^m + A_t^{j_t} x_t = b_t^{j_t}\}$
\ENDFOR

\STATE
\STATE (Backward Step)
\STATE $\mathfrak{V}_{T + 1}^{m + 1}(\cdot)\leftarrow 0$

\FOR{$t = T, \ldots, 2$}

\FOR{$j = 1, \ldots, N$}

\STATE $\pi_t^j\leftarrow \text{argmax}_\pi \min_{x_t\ge 0}\ (c_t^j)^{\sT} x_t + \mathfrak{V}_{t + 1}^{m + 1}(x_t) + \pi^{\sT}(B_t^j x_{t - 1}^m + A_t^j x_t - b_t^j)$
\STATE $\alpha_t^j \leftarrow \min_{x_t\ge 0} \{(c_t^j)^{\sT} x_t + \mathfrak{V}_{t + 1}^{m + 1}(x_t) : B_t^j x_{t - 1}^m + A_t^j x_t = b_t^j\}$

\ENDFOR

\STATE sort $\{\alpha_t^j\}_j$ so that $\alpha_t^{(j_1)}\le \alpha_t^{(j_2)}$ for $j_1\le j_2$
\STATE \small{$s_t^m(x_{t-1})\leftarrow \frac{1 - \lambda}{N} \cdot \sum_{j = 1}^N [\alpha_t^j + \ip{-(B_t^j)^\top \pi_t^j, x_{t-1} - x_{t - 1}^m}] + \lambda_t \cdot [\alpha_t^{(\kappa)} + \ip{-(B_t^{(\kappa)})^\top \pi_t^{(\kappa)}, x_{t-1} - x_{t - 1}^m}]$}
\STATE \tab \tab \small{$+ \frac{\lambda}{\alpha N} \cdot \sum_{j = \kappa + 1}^N [(\alpha_t^{(j)} - \alpha_t^{(\kappa)}) + \ip{-(B_t^{(j)})^\top \pi_t^{(j)} + (B_t^{(\kappa)})^\top \pi_t^{(\kappa)},\ x_{t-1} - x_{t - 1}^m}]$}
\STATE $\mathfrak{Q}_t^{m+1}(\cdot)\leftarrow \max\{\mathfrak{Q}_t^m(\cdot), s_t^m(\cdot)\}$

\ENDFOR

\STATE $m\leftarrow m + 1$

\ENDWHILE
\end{algorithmic}
\end{algorithm}

\subsection{Identifying ``bad'' outcomes} \label{sec:bo}
As mentioned in the last paragraph of section \ref{sec:mult}, we can find the corresponding worst case probability  density by identifying the bad outcomes. To identify the bad outcomes, we need to verify  the order of values of $V_t(\bar{x}_{t - 1}, \xi_t^j)$, $j=1,\ldots,N$, at each stage $t = 2, \ldots, T$. Although the required values   $ V_t(\bar{x}_{t - 1}, \xi_t^j)$ are not available, we have access to the approximations $\uV_t^{m+1}(x_{t-1}^m, \xi_t^j)$ at the $m$-th iteration  of the SDDP algorithm. If the approximations are good enough, by continuity argument  the order of $\uV_t^{m+1}(x_{t-1}^m, \xi_t^j)$  is more or less the same as the order of $V_t(\bar{x}_{t-1}, \xi_t^j)$, $j=1,...,N$. Since the approximations are improved in each iteration of the SDDP algorithm, we expect the order of $\uV_t^{m+1}(x_{t-1}^m, \xi_t^j)$ to stabilize after a certain number of iterations. In particular, if the SDDP algorithm is run for a sufficient number of iterations, then those bad outcomes should be the outcomes that correspond to higher values of $\uV_t^{m+1}(x_{t-1}^m, \xi_t^j)$ (in each iteration) the most frequently.

Following this idea, we try to identify the bad outcomes via frequency, i.e., by counting the number of iterations that  values
$$
\nu_{t,j}^{m+1}:=\uV_t^{m+1}(x_{t-1}^m, \xi_t^j),\;j=1,...,N,
$$
are  large at  a given  stage $t$. Formally,  let  $\nu_{t,(1)}^{m+1}\le \ldots\le \nu_{t,(N)}^{m+1}$ be these values in the increasing order,
 $\kappa := \lceil (1 - \alpha) N \rceil$ and
\begin{equation} \label{kappa}
 \nu_{t,(\kappa)}^{m+1}  =\valrisk_{\alpha}\left(\nu_{t}^{m+1}\right)
\end{equation}
 be the corresponding quantile. Let $|A|$ denote the cardinality of a set $A$, and define
\begin{equation} \label{bo-1}
W_t^j := \left |\big \{m : \uV_t^{m+1}(x_{t-1}^m, \xi_t^j)\ge  \nu_{t,(\kappa)}^{m+1}\big   \}\right |,\;  1\le t\le T,\; 1\le j\le N.
\end{equation}
In plain words, for each $t$ and $j_0$, $W_t^{j_0}$ is the frequency (i.e., the number of iterations) that $\uV_t^{m+1}(x_{t-1}^m, \xi_t^{j_0})$ is one of the $\lceil\alpha N\rceil$-largest values in $\{\uV_t^{m+1}(x_{t-1}^m, \xi_t^j)\}_{1\le j\le N}$. Ideally, the bad outcomes $\xi_t^j$ should correspond to a high value $W_t^j$. For the operation planning problem studied in section \ref{sec:ne}, numerical experiments indicate that $W_t^j$ can be as high as $2900$, after a total of $3000$ iterations,  in most of the stages.

We next discuss how to assign probability weights to outcomes $\xi_t^j$ based on values of $W_t^j$. Let   $W_t^{(1)}\le \ldots\le W_t^{(N)}$ be values $W_t^j$  arranged in the increasing order and  $q_t^j:=\cq_j(W_t)$, $j=1,\ldots,N$,  be values defined in \eqref{distr} associated with vector $W_t=(W^1_t,\ldots,W_t^N)$, then we assign weights $q_t^j$ to $\xi_t^j$.
If $W_t^{(j)}$ are  not in the  strictly increasing order (i.e., some of these values are equal to each other), we still  may   sort $W_t^j$   to ensure  that $\sum_{j = 1}^N q_t^j = 1$ for $t = 1,\ldots, T$.

\begin{remark}
\label{rem-acce}
{\rm
Suppose the frequencies  $W_t^j$ are given, and the weights $q_t^j$ are assigned to outcomes $\xi_t^j$ based on the values of $W_t^j$. Let each scenario $\{\xi_t^{j_t}\}_{t = 2, \ldots, T-1}$ be sampled at forward steps with probability $\prod_{t = 2}^{T-1} q_t^{j_t}$ (we refer to this as the {\em biased sampling}). Then the forward steps produce trial points $\{x_t^m\}_{t = 1, \ldots, T-1}$ that are different from the ones obtained under uniform sampling (as described in \eqref{approx2}).  {\em It turns out that  the convergence of the risk averse SDDP algorithm is improved  when we utilize those points at backward steps.}

The improvement in convergence can be explained in the following way.  By \eqref{avr2} we have that
\[
\begin{array}{ll}
\rho_t(\uV_t^m(x_{t-1}^m, \xi_t)) = \sum_{j = 1}^N q_{t,j}^m \nu_{t,(j)}^{m},
\end{array}
\]
with $q_{t,j}^m=\cq_j(\nu_t^m)$.  When the order   of $\{\uV_t^m(x_{t-1}, \xi_t^j)\}$ is  stable, sampling $\xi_t^j$ that corresponds to  higher weights  $q_{t,j}$ increases chances of appearing  bad scenarios. This is reminiscent of the Importance Sampling (IS) techniques, although is not exactly the same.
}
\end{remark}

\begin{remark}[Change-of-measure risk neutral problem]
\label{rem-change}
{\rm
We can construct a new risk neutral problem, referred to as the {\em change-of-measure risk neutral problem}, such that the outcomes $\{\xi_t^j\}$ at stage $t$ are  assigned the  respective weights  $\{q_t^j\}$. Such problem can be solved by applying the risk neutral SDDP algorithm with uniform sampling at forward steps. Here, uniform sampling means that  each scenario $\{\xi_t^{j_t}\}_{t = 2, \ldots, T-1}$ is sampled with probability $\prod_{t = 2}^{T-1} q_t^{j_t}$ rather than $1/N^{T-2}$, since the weights of $\xi_t^j$ are $q_t^j$ instead of $1/N$.

In section \ref{sec:ne}, we apply the risk neutral SDDP algorithm (with $\lambda = 0$) to the change-of-measure risk neutral problem constructed for the operation planning problem. The risk neutral SDDP algorithm yields similar approximations as the ones obtained by applying the risk averse SDDP algorithm to the original risk averse formulation of the operation planning problem. This offers an explanation of the nested risk averse formulation of the objective function, that such formulation is equivalent to the risk neutral problem where the risk is implicitly controlled by assigning higher weights to ``bad scenarios''.
}
\end{remark}

\subsection{Dynamic Biased Sampling} \label{sec:dbs}

In Remark \ref{rem-acce}, we have discussed how a biased sampling method for generating scenarios (thus trial points) at forward steps can help to  improve convergence of the risk averse version of the  SDDP algorithm. However, those bad outcomes are not known a priori. In order to acquire the frequencies $W_t^j$, we need to run the risk averse SDDP algorithm with uniform sampling first, which can be time consuming given the size of the problem. For the purpose of improving the convergence of the risk averse SDDP algorithm, we can instead employ a dynamic biased sampling method (for generating scenarios at forward steps) that allocates more weights to a dynamic set of outcomes. Such set gets updated in each iteration throughout the algorithm. As illustrated in section \ref{sec:ne}, the performance of the risk averse SDDP algorithm equipped with dynamic biased sampling is similar to the performance of the algorithm with biased sampling.

We next discuss how the dynamic set is updated. Consider variables $\gamma_t^j, t = 1, \ldots, T, j = 1, \ldots, N$, that represents the adjusted frequency of outcomes $\xi_t^j$. Initially, set $\gamma_t^j = 0$ for all $t, j$. At the backward step of iteration $m$, we increase $\gamma_t^j$ by $1$ if $\uV_t^{m+1}(x_{t-1}^m, \xi_t^j)\ge \nu_{t,(\kappa)}^{m+1}$. We then set $\gamma_t^j\leftarrow \gamma_t^j \cdot \frac{m}{m + 1}$ for all $t$ and $j$, and proceed to the next iteration. The adjustment weakens the impact of early iterations, where the approximations are still inaccurate. The dynamic set of bad outcomes in each stage $t$ are the outcomes with the highest adjusted frequency so far, and weights are assigned according to  \eqref{distr}. When the bad outcomes are  stable, the dynamic set should coincide with the set of those outcomes.
A complete description of the algorithm can be found in Appendix A.

\section{Numerical Experiments} \label{sec:ne}

The numerical experiments are performed on an aggregated representation of the Brazilian Interconnected Power System operation planning problem  with historical data as of January 2011. The study horizon is  60 stages and the total number of considered stages  is $T = 120$. The scenario tree has 100 realizations at  every stage (i.e., $N = 100$), and each realization $\xi_t^j, t = 2, \ldots, 120, j = 1, \ldots, 100$, has the same weight $1/N = 1/100$.  The random data process is represented by  four dimensional vectors of monthly inflows, aggregated in four regions,  and modeled as the first order periodical time series process.   The total number of state variables is 8.
We refer to \cite{ejor}  for the detail description of the model (see also Appendix B).
We implement the risk averse SDDP algorithms with different sampling methods (for generating scenarios and  thus trial points). We also implement the risk neutral SDDP algorithm for the constructed change-of-measure risk neutral problem (see Remark \ref{rem-change}). Both implementations were written in C++ and using Gurobi 8.1. Dual simplex was used as a default method for the LP solver.


We conduct the numerical experiments with parameters $\lambda = 0.2$ and $\alpha = 0.05$ (see \eqref{cm-crm}) in the following three steps:

\begin{enumerate}

\item The first step is to run the risk averse SDDP algorithm with uniform sampling, referred to as  ``raus" (where each scenario $\{\xi_t^{j_t}\}_{t = 2, \ldots, T-1}$ is sampled with probability $1/N^{T-2}$  at forward steps) and to identify the bad outcomes as discussed in subsection \ref{sec:bo}.

\item The second step is to run the risk averse SDDP algorithm with biased sampling, referred to as  ``rabs" (where each scenario $\{\xi_t^{j_t}\}_{t = 2, \ldots, T-1}$ is sampled with probability $\prod_{t = 2}^{T-1} q_t^{j_t}$, see   Remark \ref{rem-acce}) and with dynamic biased sampling ``radbs" (see section \ref{sec:dbs}). We then compare lower bounds of the optimal objective value of the risk averse problem produced by the SDDP algorithm with different sampling methods (recall that  the lower bound is given by the quantity $f_1(x_1^m) + \mathfrak{V}_2^m(x_1^m)$).

\item \label{step3} The final step is to solve the change-of-measure risk neutral problem by the risk neutral SDDP algorithm with uniform sampling, referred to as  ``nrn" (where each scenario $\{\xi_t^{j_t}\}_{t = 2, \ldots, T-1}$ is sampled with probability $\prod_{t = 2}^{T-1} q_t^{j_t}$, see  Remark \ref{rem-change}). In addition, we compare the policies corresponding to the risk averse problem and the change-of-measure risk neutral problem.

\end{enumerate}

{\color{black} In addition, we perform numerical experiments in the extreme case where $\lambda = 0.5$ and $\alpha = 0.05$} to examine the effectiveness of dynamic biased sampling in improving the convergence of the lower bounds.

\paragraph{Experiment with $\lambda = 0.2, \alpha = 0.05$.}


We first run the risk averse SDDP algorithm with uniform sampling for 3000 iterations. As discussed in subsection \ref{sec:bo}, the bad outcomes are identified via frequencies $W_t^j$, and the weights $q_t^j$ are assigned accordingly. Here, $\kappa = \lceil (1 - \alpha) N \rceil = 95$, and more weights are assigned to the outcomes corresponding to the top five frequencies.

\begin{figure}[H]
  \centering
    \includegraphics[scale=0.27]{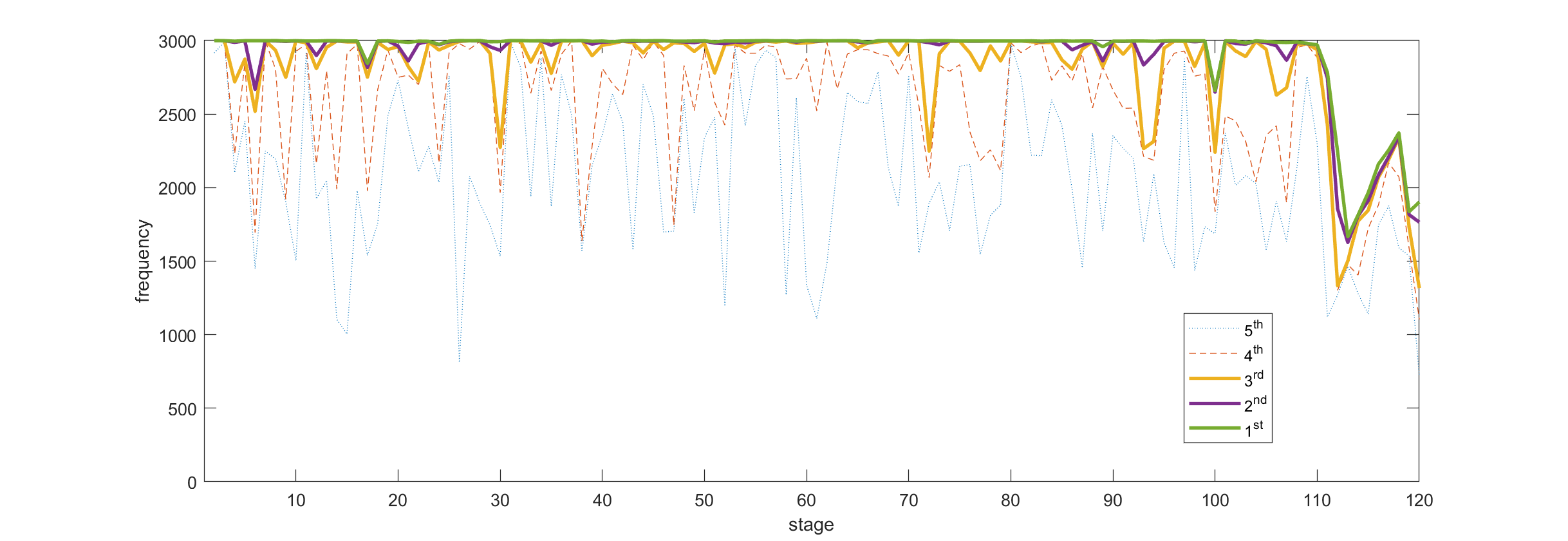}
	\caption{The counts  of the top five frequencies  at each stage after 3000 iterations.}
	\label{fig1}
\end{figure}

Figure \ref{fig1} plots the top five frequencies at each stage after 3000 iterations. In particular, the $i$-th curve corresponds to the $i$-th highest frequency by stage, which are the points $\{(t, W_t^{(N - i + 1)})\}_{t = 2, \ldots, T}$ in the notation in subsection \ref{sec:bo}. A high frequency suggests the corresponding outcome is indeed a bad outcome and should be assigned more weight in the worst case probability density. Observe that the frequency of the top three curves are almost the same as the total number of iterations. The fluctuations of frequency of the $4$ and $5$-th curve were partly due to inaccuracy of the approximations in the early iterations, and the corresponding outcomes become stable as the number of iterations increases. Such evidence supports the validity of the change-of-measure approach.

\begin{figure}[H]
  \centering
    \includegraphics[scale=0.35]{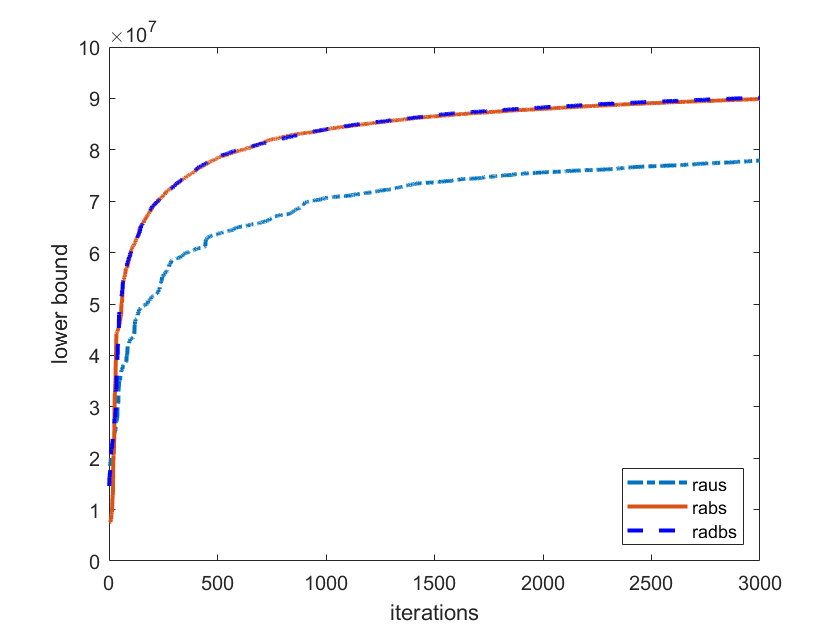}
	\includegraphics[scale=0.35]{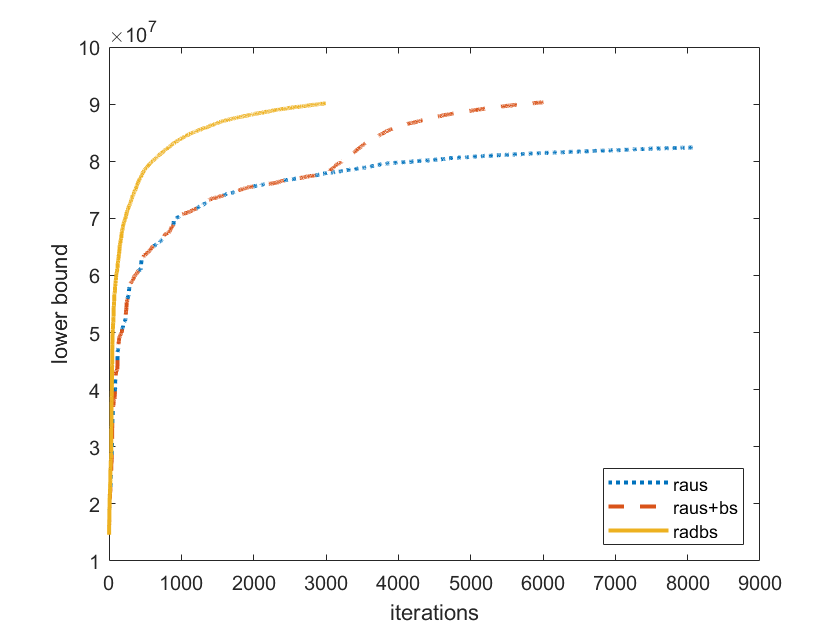}
	\caption{convergence of lower bounds of objective values produced by SDDP algorithm with uniform and (dynamic) biased sampling methods ($\lambda = 0.2$ and $\alpha = 0.05$).}
	\label{fig2}
\end{figure}

Figure \ref{fig2} contains two plots of lower bounds of the optimal objective value of the considered  risk averse problem produced by the risk averse SDDP algorithm with different sampling methods (for generating scenarios thus trial points at forward steps), i.e., for each sampling method, we plot the points $(m, f_1(x_1^m) + \mathfrak{V}_2^m(x_1^m))$ by iteration $m$ (see also section \ref{sec:sddp}).

The plot on the left of Figure \ref{fig2} compares ``raus", ``rabs", and ``radbs'', which correspond to the risk averse  SDDP algorithm where the trial points are generated under uniform, biased, and dynamic biased sampling, respectively. It shows the convergence of ``rabs'' and ``radbs'' are almost the same, which are much better than that of ``raus''. As discussed in subsection \ref{sec:dbs}, ``rabs'' requires weights obtained by first running ``raus'', which can be time consuming given the size of the problem. This makes ``radbs'' a better choice for the purpose of improving the convergence of the risk averse SDDP algorithm. Nevertheless, the approximations generated by ``rabs'' and ``radbs'' are similar, and they are better than the one generated by ``raus''. In the remaining section, we use the approximation $\{\mathfrak{V}_t^A\}$ generated  by ``radbs'' after 3000 iterations as the reference approximation for the risk averse problem.

The plot on the right of Figure \ref{fig2} compares ``raus", ``radbs ", and ``raus+bs" in varying iterations, where ``raus+bs'' is the risk averse SDDP algorithm equipped with uniform sampling in the first 3000 iterations and switched to biased sampling afterwards. We see   that the lower bound generated by ``raus'' in 8000 iterations is still much worse than the one generated by ``rabs'' in 3000 iterations, and the gap closes slowly. Besides, ``raus+bs'' takes about 6000 iterations to reach what ``radbs'' achieves in 3000 iterations.

As the final step, we solve the change-of-measure risk neutral problem by ``nrn" (see step \ref{step3}). Let $\{\mathfrak{V}_t^N\}$ denote the approximations produced by ``nrn" after 3000 iterations. We next compare the policies generated by the approximations $\{\mathfrak{V}_t^A\}$ and $\{\mathfrak{V}_t^N\}$. Formally, for a given scenario $\{\xi_t^{j_t}\}_{t = 1, \ldots, T}$ and approximations $\{\mathfrak{V}_t\}_{t = 1, \ldots, T}$, a policy generated by the approximations is $\{x_t(\xi_t^{j_t})\}_{t = 1, \ldots, T}$ such that
\begin{equation} \label{policy}
x_t(\xi_t^{j_t}) \in \arg\min\limits_{x_t\in \mathcal{X}_t(x_{t-1}^m, \xi_t^{j_t})} \{f_t(x_t, \xi_t^{j_t}) +  \mathfrak{V}_{t+1}(x_t)\}.
\end{equation}
The policies are the decisions to be implemented given the scenario. If $\{\mathfrak{V}_t^A\}$ and $\{\mathfrak{V}_t^N\}$ yield similar policies for each scenario, then the risk averse formulation is similar to the change-of-measure risk neutral formulation, at least from a numerical perspective. Instead of comparing policies for all scenarios, we randomly sample 3000 scenarios and compare the policies by plotting their paths (or distribution).

\begin{figure}[H]
  \centering
	\includegraphics[scale=0.5]{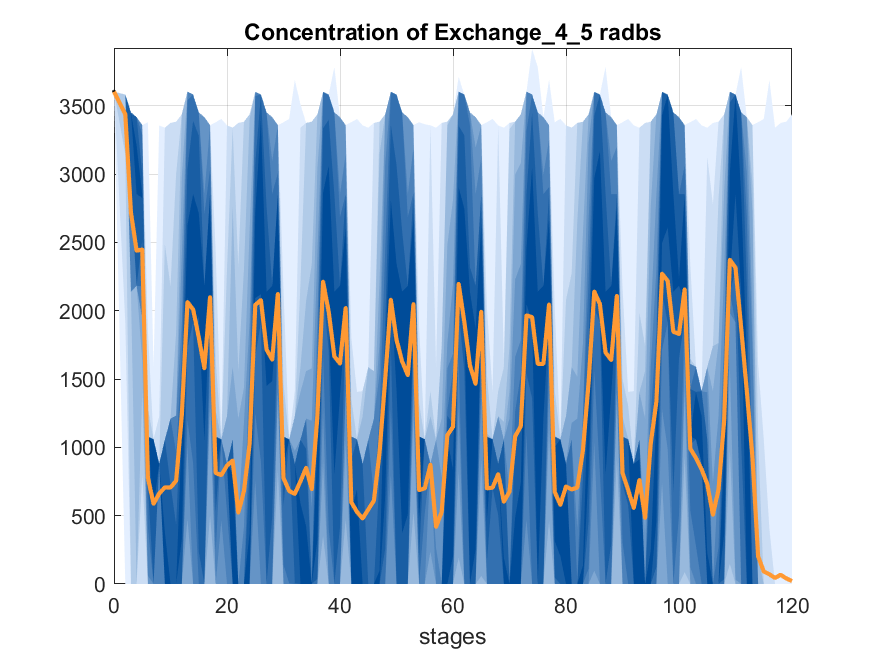}
    \includegraphics[scale=0.5]{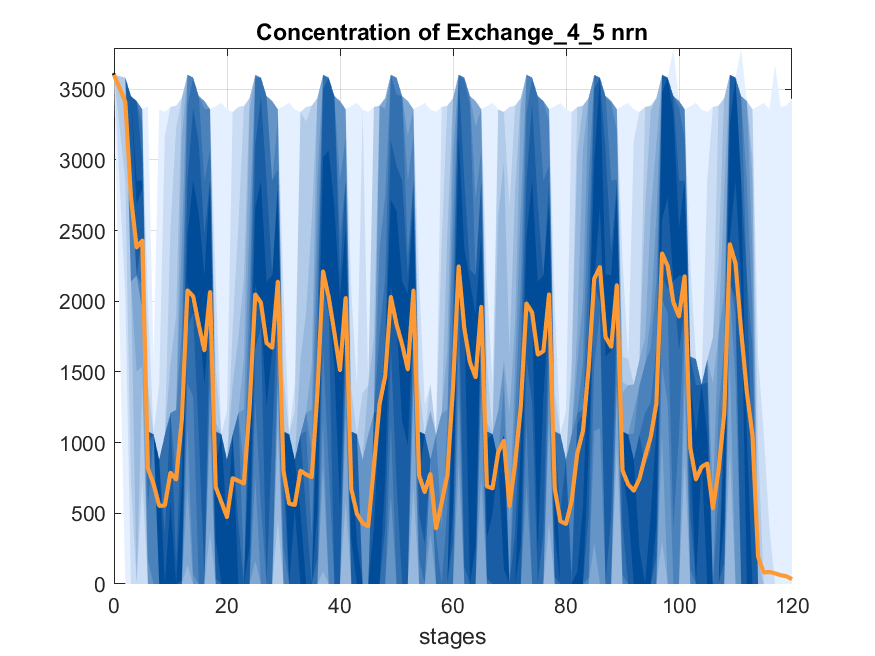}

	\caption{fanplot of paths of stage variable 1 generated by $\{\mathfrak{V}_t^A\}$ and $\{\mathfrak{V}_t^N\}$.}
	\label{fig3}
\end{figure}

\begin{figure}[H]
  \centering
    \includegraphics[scale=0.5]{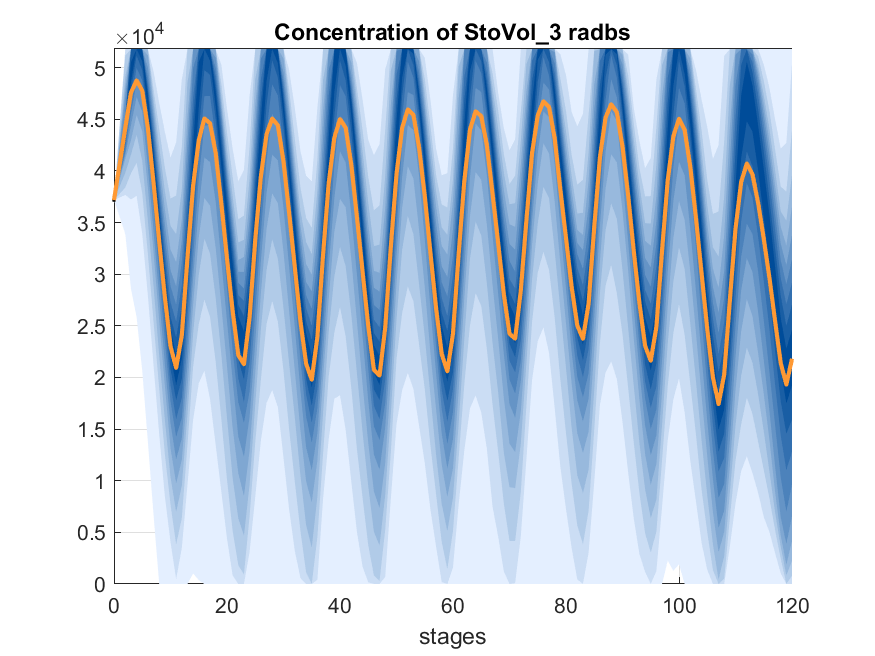}
    \includegraphics[scale=0.5]{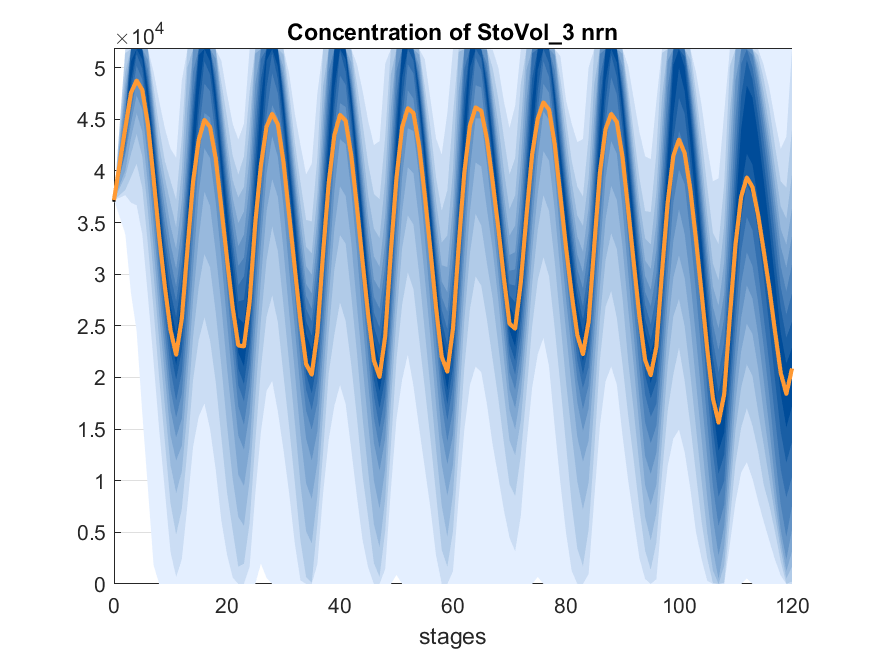}

	\caption{fanplot of paths of stage variable 2 generated by $\{\mathfrak{V}_t^A\}$ and $\{\mathfrak{V}_t^N\}$.}
	\label{fig4}
\end{figure}

\begin{figure}[H]
  \centering
    \includegraphics[scale=0.5]{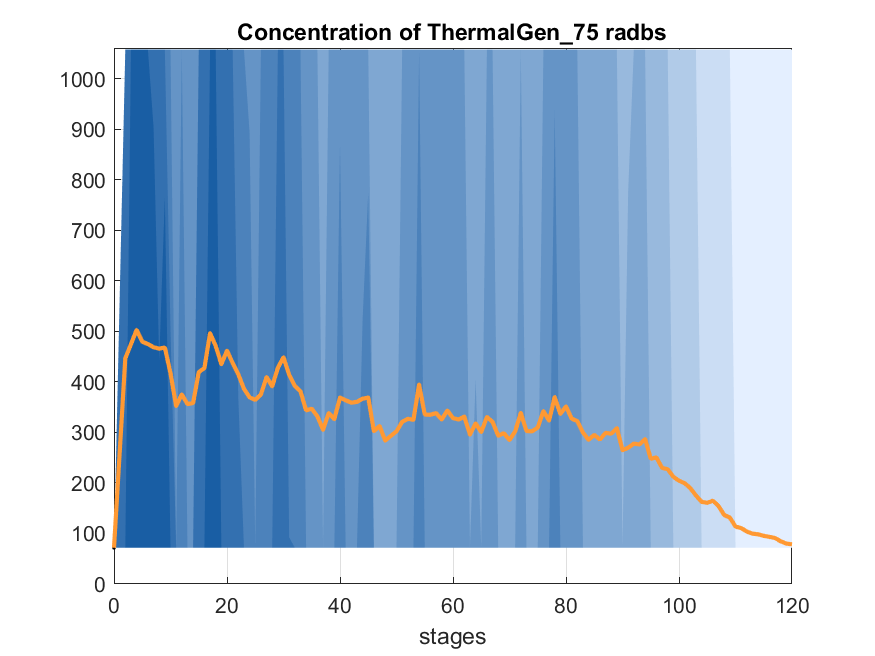}
    \includegraphics[scale=0.5]{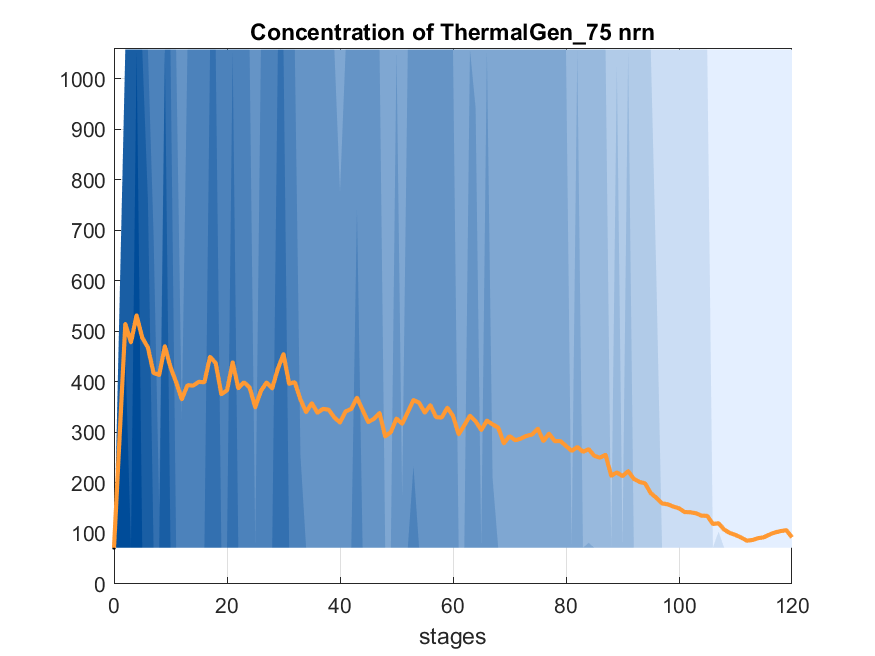}

	\caption{fanplot of paths of stage variable 3 generated by $\{\mathfrak{V}_t^A\}$ and $\{\mathfrak{V}_t^N\}$.}
	\label{fig5}
\end{figure}

For operation planning problem, the set of stage variables is  the same across all stages, and each stage variable at stage $t$ corresponds to an entry of $x_t(\xi_t^{j_t})$. We sample 3000 scenarios uniformly at random and generate a fanplot for paths of each variable and of stage objective value (explained below). In Figures  \ref{fig3}, \ref{fig4} and \ref{fig5}, we present fanplots that exhibit typical behavior of the variables. The dark area on the fanplots is where the paths of the variables are highly concentrated, whereas the light area is the opposite. The orange curve on the plots is the average of the variables by stage.

The plots on the left of Figures \ref{fig3}, \ref{fig4} and \ref{fig5} correspond to the risk averse problem, whereas the ones on the right correspond to the constructed risk neutral problem. We can see from the plots that concentrations of policies of the same variable generated by $\{\mathfrak{V}_t^A\}$ and $\{\mathfrak{V}_t^N\}$ are similar.

\begin{figure}[H]
  \centering
    \includegraphics[scale=0.5]{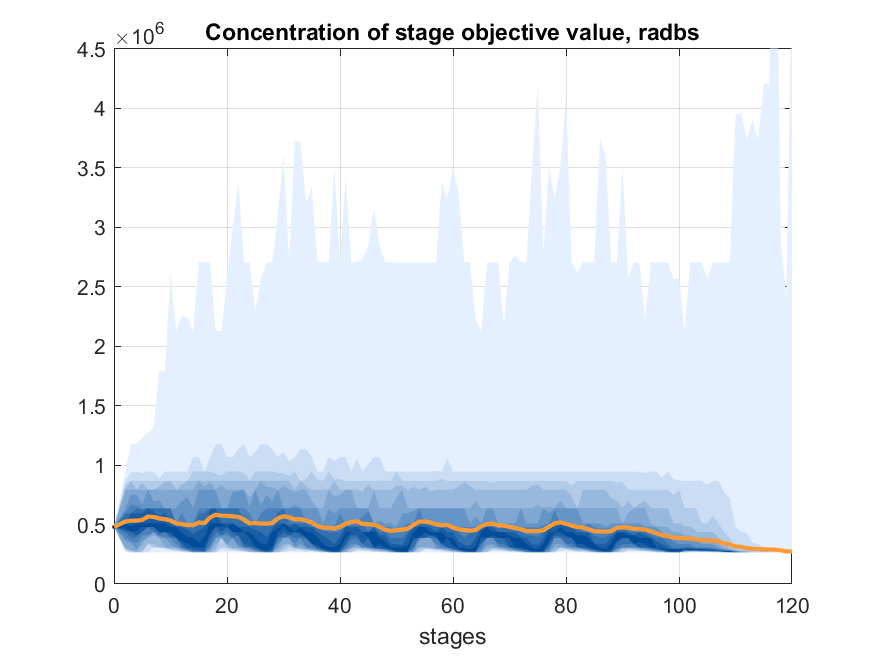}
    \includegraphics[scale=0.5]{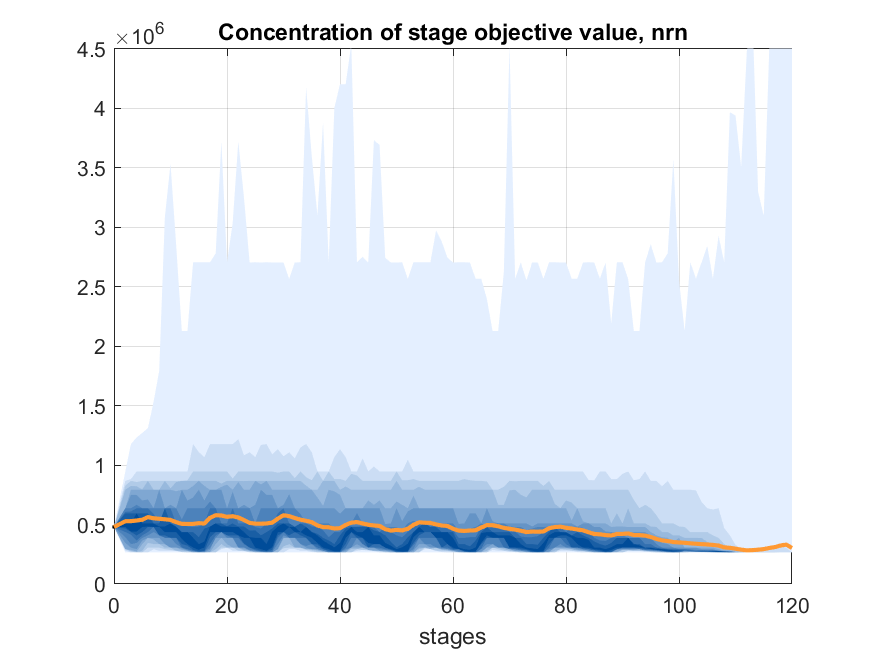}

	\caption{fanplot of paths of stage objective values generated by $\{\mathfrak{V}_t^A\}$ and $\{\mathfrak{V}_t^N\}$.}
	\label{fig6}
\end{figure}

Figure \ref{fig6} plots concentrations of the stage objective value $f_t(x_t(\xi_t^{j_t}), \xi_t^{j_t})$ for policies $x_t(\xi_t^{j_t})$ generated by $\{\mathfrak{V}_t^A\}$ and $\{\mathfrak{V}_t^N\}$, respectively. Figures \ref{fig3} to \ref{fig6} indicate that the risk averse and the constructed risk neutral problem yield similar policies and objective values, hence the risk averse formulation and the change-of-measure risk neutral formulation are similar.

\begin{figure}[H]
  \centering
    \includegraphics[scale=0.27]{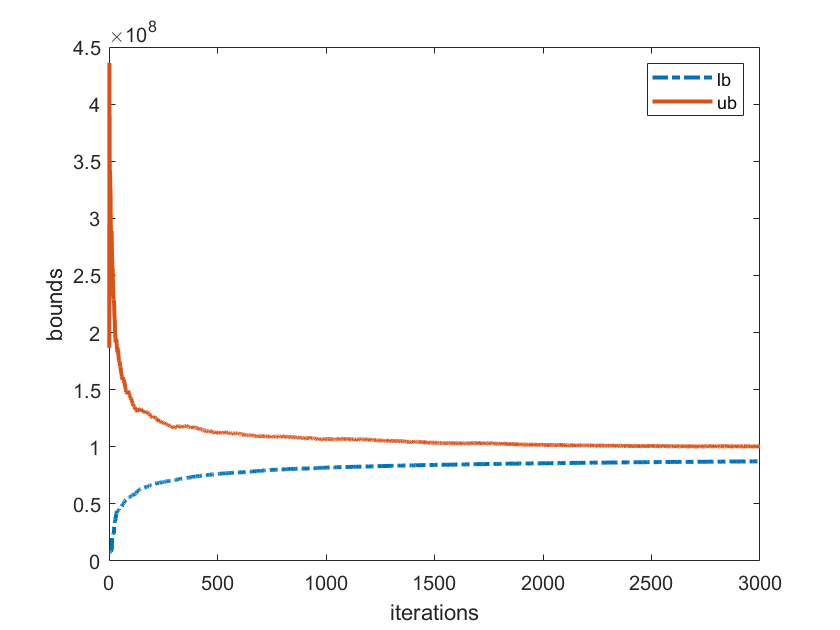}

	\caption{lower bounds and statistical upper bounds per iteration  of the optimal  objective value  of the change-of-measure risk neutral problem.}
	\label{fig7}
\end{figure}

Figure \ref{fig7} plots the lower bounds and statistical upper bounds of the optimal objective value of the change-of-measure risk neutral problem. Here, the statistical upper bound at iteration $m$ is defined to be the quantity $\frac{1}{m} \sum_{i = 1}^m \sum_{t = 1}^T f_t(x_t^i (\xi_t^{j_t^i}), \xi_t^{j_t^i})$, where $\xi_t^{j_t^i}$, $t=1,\ldots,T$,  are  the scenarios  generated at iteration $i$ (note that it is different from the statistical upper bound introduced in \cite{per91}). After 3000 iterations the gap between the lower and upper bounds is about 14\%.
One advantage of applying the SDDP algorithm to the constructed  risk neutral problem is the possibility of incorporating the termination criterion based on the  lower  and   upper bounds.

\paragraph{Experiment with $\lambda = 0.5, \alpha = 0.05$.} {\color{black} The case $\lambda = 0.5$ leads to an extreme risk aversion}. The parameter $\lambda = 0.5$ was chosen instead of some larger values (e.g. $\lambda = 1$) to better reflect the reality.

\begin{figure}[H]
  \centering
    \includegraphics[scale=0.35]{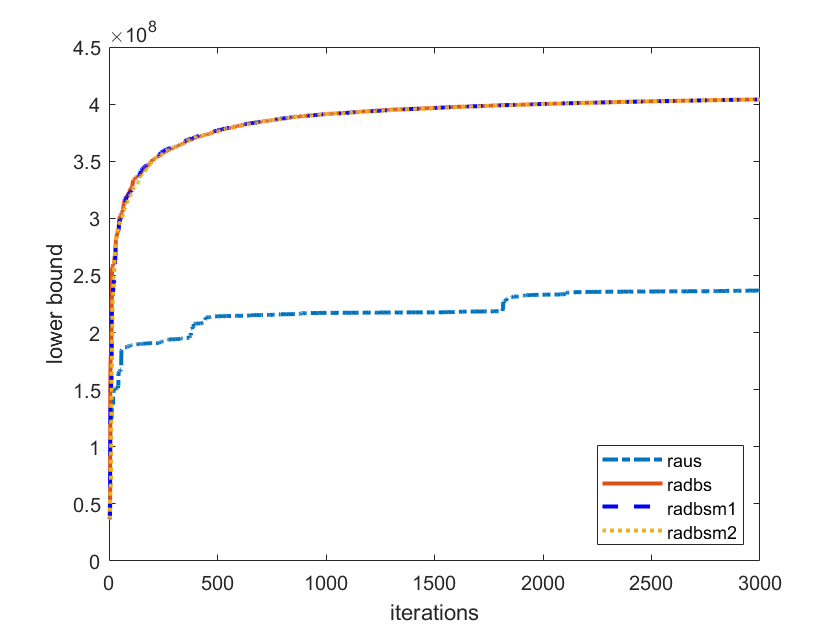}

	\caption{convergence of lower bounds of objective values produced by SDDP algorithm with uniform and various dynamic biased sampling methods ($\lambda = 0.5$ and $\alpha = 0.05$).}
	\label{fig8}
\end{figure}

Figure \ref{fig8} compares the convergence of lower bounds generated by ``raus'', ``radbs'', and two variants of ``radbs'', namely ``radbsm1'' and ``radbsm2''. Those two variants differ from ``radbs'' in the adjustment $\gamma_t^j \leftarrow \gamma_t^j \cdot \frac{m}{m + 1}$ in iteration $m$; ``radbsm1'' simply does not have such adjustment, whereas the adjustment in ``radbsm2'' is replaced by $\gamma_t^j \leftarrow \gamma_t^j \cdot (1 - 0.5^m)$. From Figure \ref{fig8}, we can see the convergence of the three dynamic biased sampling methods is more or less the same. However, they achieve a striking speed up over ``raus'', where their lower bounds are almost twice the lower bound of ``raus'' after 3000 iterations. Hence, the dynamic biased sampling method gets more effective when $\lambda$ increases. {\color{black} Note in the other extreme case where $\lambda = 0$, biased sampling reduces to uniform sampling, thus the dynamic biased sampling approach does not offer any speed up.}

\section{Conclusions} \label{sec:conc}


We discuss  the risk averse (multistage) stochastic programming and  the idea of ``bad'' outcomes. We show   how to identify  the ``bad'' outcomes  and use this with the SDDP method.
Numerical experiments were conducted on the Brazilian Interconnected Power System operation planning problem, with  the results  presented and discussed in section \ref{sec:ne}. The convergence of lower bounds generated by the risk averse version of the  SDDP algorithm with different sampling methods  was examined.  It was observed  that
the  (dynamic) biased sampling approach for generating trial points considerably improved performance of the lower bounds.
 We also compared the
solutions generated by the risk averse version of the SDDP method and the changed-of-measure  risk neutral  SDDP method
and concluded that the policies and the  objective values generated by these  two approaches  are similar.

\paragraph{Acknowledgement}
The authors are indebted to Filipe Goulart Cabral and anonymous referees   for constructive comments which helped to improve presentation of the manuscript.

\bibliographystyle{plain}

\bibliography{references}

\section{Appendix A} \label{appen:A}

{\setstretch{1.0}

\begin{algorithm}[H]
\caption{the SDDP algorithm for linear problem (dynamic biased sampling)}
\begin{algorithmic}
\STATE $\kappa\leftarrow \lceil(1 - \alpha)N)\rceil, m\leftarrow 1, \mathfrak{V}_t^1(\cdot)\leftarrow 0\ \forall t, \gamma_t^j\leftarrow 0\ \forall t, j$

\WHILE{termination criterion not met}

\IF{$m = 1$}
\STATE sample $\{\xi_t^{j_t}\}_{t = 2, \ldots, T-1}$ with probability $N^{2-T}$
\ELSE
\FOR{$t = 2, \ldots, T$}
\STATE sort $\{\gamma_t^j\}_j$ so that $\gamma_t^{(j_1)}\le \gamma_t^{(j_2)}$ for $j_1\le j_2$
\FOR{$j = 1, \ldots, N$}
\STATE $w_t^j \leftarrow \begin{cases} (1 - \lambda)/N & \text{if } \gamma_t^j < \gamma_t^{(\kappa)} \\ (1 - \lambda)/N + \lambda - \lambda(N - \kappa)/(\alpha N) & \text{if } \gamma_t^j = \gamma_t^{(\kappa)} \\ (1 - \lambda)/N + \lambda/(\alpha N) & \text{if } \gamma_t^j > \gamma_t^{(\kappa)} \end{cases}$
\ENDFOR

\STATE sample $\{\xi_t^{j_t}\}_{t = 2, \ldots, T-1}$ with probability $\prod_{t = 2}^{T-1} w_t^{j_t}$

\ENDFOR
\ENDIF

\STATE $x_1^m\leftarrow \text{argmin}_{x_1\ge 0} \{c_1^{\sT} x_1 + \mathfrak{V}_2^m(x_1) : A_1 x_1 = b_1\}$
\FOR{$t = 2, \ldots, T - 1$}
\STATE $x_t^m \leftarrow \text{argmin}_{x_t\ge 0} \{(c_t^{j_t})^{\sT} x_t + \mathfrak{V}_{t + 1}^m(x_t) : B_t^{j_t} x_{t - 1}^m + A_t^{j_t} x_t = b_t^{j_t}\}$
\ENDFOR

\STATE

\STATE $\mathfrak{V}_{T + 1}^{m + 1}(\cdot)\leftarrow 0$
\FOR{$t = T, \ldots, 2$}
\FOR{$j = 1, \ldots, N$}

\STATE $\pi_t^j\leftarrow \text{argmax}_\pi \min_{x_t\ge 0}\ (c_t^j)^{\sT} x_t + \mathfrak{V}_{t + 1}^{m + 1}(x_t) + \pi^{\sT}(B_t^j x_{t - 1}^m + A_t^j x_t - b_t^j)$
\STATE $\alpha_t^j \leftarrow \min_{x_t\ge 0} \{(c_t^j)^{\sT} x_t + \mathfrak{V}_{t + 1}^{m + 1}(x_t) : B_t^j x_{t - 1}^m + A_t^j x_t = b_t^j\}$

\ENDFOR

\STATE sort $\{\alpha_t^j\}_j$ so that $\alpha_t^{(j_1)}\le \alpha_t^{(j_2)}$ for $j_1\le j_2$
\FOR{$j = 1, \ldots, N$}
\IF{$\alpha_t^j\ge \alpha_t^{(\kappa)}$}

\STATE $\gamma_t^j\leftarrow \gamma_t^j + 1$

\ENDIF
\STATE $\gamma_t^j\leftarrow \gamma_t^j\cdot \frac{m}{m + 1}$
\ENDFOR

\STATE \small{$s_t^m(x_{t-1})\leftarrow \frac{1 - \lambda}{N} \cdot \sum_{j = 1}^N [\alpha_t^j + \ip{-(B_t^j)^\top \pi_t^j, x_{t-1} - x_{t - 1}^m}] + \lambda_t \cdot [\alpha_t^{(\kappa)} + \ip{-(B_t^{(\kappa)})^\top \pi_t^{(\kappa)}, x_{t-1} - x_{t - 1}^m}]$}
\STATE \tab \tab \small{$+ \frac{\lambda}{\alpha N} \cdot \sum_{j = \kappa + 1}^N [(\alpha_t^{(j)} - \alpha_t^{(\kappa)}) + \ip{-(B_t^{(j)})^\top \pi_t^{(j)} + (B_t^{(\kappa)})^\top \pi_t^{(\kappa)},\ x_{t-1} - x_{t - 1}^m}]$}
\STATE $\mathfrak{Q}_t^{m+1}(\cdot)\leftarrow \max\{\mathfrak{Q}_t^m(\cdot), s_t^m(\cdot)\}$

\ENDFOR

\STATE $m\leftarrow m + 1$

\ENDWHILE
\end{algorithmic}
\end{algorithm}

}


\section{Appendix B} \label{appen:B}

{\bf Long Term Operation Planning Problem}

The dynamic programming equation for the long term operation planning problem can be written as
\begin{align*}
Q_t([v_t, a_{[t-p, t-1]}], \eta_t) = \min\ \ & \sum_{k\in K}(\sum_{j\in T_k} c_j g_{tj} + \sum_{i\in U_k} \tilde{c}_i Def_{ti}) + \beta \mathcal{Q}_{t+1} [v_{t + 1}, a_{[t-p+1, t]}]) \\
			\text{s.t.}\ \ & v_{t + 1} = v_t + a_t - q_t - s_t \\
							& a_t = \text{diag}(\eta_t)(\Phi_{t, 0} + \sum_{\nu = 1}^p \Phi_{t, \nu} a_{t - \nu}) \\
							& q_{tk} + \sum_{j\in T_k} g_{tj} + \sum_{i\in U_k} Def_{ti} + \sum_{l\in \Omega_k} (f_{tlk} - f_{tkl}) = d_{tk},\tab \forall k\in K \\
							& 0\le v_{t + 1}\le \overline{v}, 0\le q_t\le \overline{q}, 0\le s_t, \\
							& \underline{g}\le g_t\le \overline{g}, 0\le Def_t\le \overline{Def}, \underline{f}\le f_t\le \overline{f},
\end{align*}
for $t = 1\ldots, T$, where
\[ \mathcal{Q}_{t+1}([v_{t+1}, a_{[t-p+1, t]}]) = \begin{cases} \rho_{t + 1} [Q_{t + 1}([v_{t+1}, a_{[t-p+1, t]}], \eta_{t+1})] & t = 1, \ldots, T - 1 \\
		0 & t = T \end{cases}, \]
and
\[ \rho_t[Z] := (1 - \lambda)\cdot \mathbb{E}[Z | \eta_{t - 1}] + \lambda\cdot AV@R_{\alpha} [Z | \eta_{t - 1}] \]
with $\lambda\in [0, 1]$ and $\alpha\in [0, 1]$ being chosen parameters. \\

The {\it objective function}
\[ \sum_{k\in K}(\sum_{j\in T_k} c_j g_{tj} + \sum_{i\in U_k} \tilde{c}_i Def_{ti}) + \beta \mathcal{Q}_{t+1} [v_{t + 1}, a_{[t-p+1, t]}]) \]
represents the sum of the total cost for thermal generation and deficit with $\mathcal{Q}_{t+1} ([v_{t + 1}, a_{[t-p+1, t]}])$, where
\begin{itemize} [noitemsep,topsep=0pt]
\item $\beta$ is a discount factor;

\item $K$ is a subsystem set;

\item $T_k$ is the thermal set for subsystem $k$;

\item $U_k$ is the deficit set for subsystem $k$. \ \\
\end{itemize}

The {\it energy balance} equation for each reservoir $k$ is
\[ v_{t+1} = v_t + a_t + q_t + s_t, \]
where
\begin{itemize} [noitemsep,topsep=0pt]
\item $v_t$ is the stored energy in the reservoir at the beginning of stage $t$;

\item $a_t$ is the energy inflow during stage $t$;

\item $q_t$ is the generated energy during stage $t$;

\item $s_t$ is the spilled energy during stage $t$. \ \\
\end{itemize}

The {\it time-series} model for the energy inflow is
\[ a_t = \text{diag}(\eta_t) (\Phi_{t, 0} + \sum_{\nu = 1}^p \Phi_{t, \nu} a_{t - \nu}), \]
where
\begin{itemize} [noitemsep,topsep=0pt]
\item $a_{t-\nu}$ is the energy inflow during stage $t - \nu$;

\item $\Phi_{t, \nu}$ is the coefficient of PVAR vector time-series model $t$;

\item $\eta_t$ is the multiplicative noise of PVAR, which independent for each stage $t$. \ \\
\end{itemize}

The {\it load balance} equation, in MW month, for each subsystem $k$ and stage $t$ is
\[ q_{tk} + \sum_{j\in T_k} g_{tj} + \sum_{i\in U_k} Def_{ti} + \sum_{l\in \Omega_k} (f_{tlk} - f_{tkl}) = d_{tk}, \]
where
\begin{itemize} [noitemsep,topsep=0pt]
\item $d_{tk}$ is load;

\item $q_{tk}$ is hydro generation;

\item $\sum_{j\in T_k} g_{tj}$ is thermal generation;

\item $\sum_{i\in U_k} Def_{ti}$ is deficit generation;

\item $\sum_{\l\in \Omega_k} (f_{tlk} - f_{tkl})$ is net energy interchange;

\item $f_{tlk}$ is the energy flow from subsystem $l$ to subsystem $k$;

\item $\Omega_k$ is the subsystems directly connected to subsystem $k$. \ \\
\end{itemize}

The bounds on variables are
\begin{itemize} [noitemsep,topsep=0pt]
\item $0\le v_{t+1}\le \overline{v}$ is the lower and upper bounds on stored energy;

\item $0\le q_t\le \overline{q}$ is the lower and upper bounds on generated energy;

\item $0\le s_t$ is non-negativity constraint of spilled energy;

\item $\underline{g}\le g_t\le \overline{g}$ is the lower and upper bounds on thermal generation;

\item $0\le Def_t\le \overline{Def}$ is the lower and upper bounds on energy deficit;

\item $\underline{f}\le f_t\le \overline{f}$ is the lower and upper bounds on energy flow. \ \\
\end{itemize}

The main idea of the deficit is to penalize the load cut by a convex piecewise linear cost function which is dependent on the load cut depth. Regarding this approach, it is important to emphasize that the deficit variable with highest associated cost is unbounded above.

For each stage $t$ the decision vector is $x_t = (v_{t + 1}, q_t, s_t, g_t, Def_t, f_t, a_t)$. In the long term operation planning problem the only considered uncertainty is the independent multiplicative noise, that is, $\xi = \eta_t$, and the cost-to-go function depends only on $[v_t, a_{[t-p, t-1]}]$ of last $p$ previous decision $x_{[t-p, t-1]}$. In this way, it is usual to write $Q_t([v_t, a_{[t-p, t-1]}], \eta_t)$ instead of $Q_t(x_{[t-p, t-1]}, \eta_t)$.




\end{document}